\documentclass[letterpaper,11pt]{article}
\usepackage{amsmath}
\usepackage{amssymb}
\usepackage{amsthm}
\newenvironment{rcase}
{\left. \begin{aligned}}
{\end{aligned}\right\rbrace}
\author{Svetlana Ekisheva\thanks{Schools of Biology and $^\dagger$Mathematics,
Georgia Institute of Technology, 
Atlanta, GA 30332-0230 USA
sveta.ekisheva@bme.gatech.edu;
houdre@math.gatech.edu} \ and Christian Houdr\'{e}$^\dagger$}

\title{Transportation Distance and the Central Limit Theorem}

\begin{document}
\maketitle
\begin{abstract}
For probability measures
on a complete separable metric space, we present sufficient conditions
for the existence of a solution
to the Kantorovich transportation problem. 
We also obtain sufficient conditions (which sometimes also become necessary)  for the convergence,
in transportation, of probability 
measures when the cost function is  continuous, non-decreasing and depends
on the distance.
As an application, the CLT in the transportation distance is proved
for independent and some dependent stationary sequences. 
\end{abstract}


\section{Introduction}

Let $(M,d)$ be a metric space and let $c:M\times M\rightarrow {\bf R}$,
be a non-negative Borel  function.
The transportation $c$-distance $T_{c}(\mu , \nu )$ between two probability measures
$\mu $ and $\nu $ defined on the Borel $\sigma $-field $\mathcal B(M)$ is given via
$$T_{c}(\mu ,\nu )=\inf {\bf E} c(X,Y).$$
Above, the infimum is taken over all $M$-valued random elements $X$ and $Y$ defined on the
probability space $(\Omega, \mathcal F, \bf{P})$
and having, respectively, $\mu $ and $\nu $ for probability distribution.
In  other words,
\begin{equation}
T_{c}(\mu ,\nu )=\inf _\Pi \int c(x,y)d\pi (x,y),
\end{equation}
where the infimum is taken over the set $\Pi$ of all probability measures on $\mathcal B(M \times M)$
with marginals $\mu$ and $\nu$. The transportation
distance is related to the celebrated Kantorovich transportation problem: if $\mu $ and $\nu $ are two
distributions of mass and if $c(x,y)$ represents
the cost of transporting a unit of mass from the location $x$ to 
the location $y$, what is the minimal total transportation cost to transfer $\mu $ to $\nu $?
The minimal total transportation cost is exactly the transportation 
distance corresponding to the cost function $c$.

The $c$-transportation distance with $c(x,y)=d^p(x,y)$,
$p\geq 1$, is associated to the Wasserstein or  Mallows $p$-distance $W_p$,
$W_p(\mu, \nu)=(T_{d^p}(\mu ,\nu ))^{1/p}$.
If $M$
is the real line ${\bf R}$ with the Euclidean distance,
the Wasserstein-Mallows $p$-distance between
two distribution functions $F$ and $G$ has the following useful representation
\begin{equation}
W^p_p(F,G)=\int _0^1 |F^{-1}(t)-G^{-1}(t)|^p dt,
\end{equation}
where the inverse transformation of $F$ is defined as
$$F^{-1}(t)=\sup \{x\in {\bf R}:F(x)\leq t\}.$$   
The representation (2) was obtained when $p=1$ by Salvemini \cite{Sal} (for discrete distributions) and by Dall'Aglio \cite{D56}
(in the general case), while for $p=2$ it is due to Mallows \cite{Mal}. It implies that the random
variables $X=F^{-1}(U )$ and $Y=G^{-1}(U)$, where  $U$ is  a uniform random variable on $(0,1)$,
are minimizers of the total
transportation cost in the transportation problem.
Major \cite{Maj} generalized (2) to a convex cost function
$c(x,y)=c(x-y)$:
\begin{equation*}
T_c(F,G)=\int _0^1 c(F^{-1}(t)-G^{-1}(t))dt.
\end{equation*}
The representation (2) is an important tool in proving
the following convergence result. 

Let $p=1,2$ and
let $F_n$, $F$ be distribution functions on ${\bf R}$
such that for any $n$, 
$\int |x|^pdF_n<+\infty$, and $\int |x|^pdF<+\infty$. Then
\begin{equation}
W_p(F_n,F)\rightarrow 0 \Longleftrightarrow
 \begin{cases}
   (a)\  F_n\Longrightarrow F,\\
   (b) \int |x|^pdF_n \rightarrow \int |x|^p dF.
   \end{cases}
\end{equation}
For $p=1$ the equivalence (3) was proved by Dobrushin \cite{Dob}, while
for $p=2$ it is due to Mallows \cite{Mal}.

\smallskip
Bickel and Freedman \cite{BF} extended the statement (3) to probability
measures $\mu _n $ and $\mu $ defined on a separable Banach space $({\bf B}, \|\cdot\|)$
and to all $p \in [1,+\infty)$ as follows:

 Let $1 \leq p < \infty$, and let $\int \|x\|^p \mu _n (dx) < \infty$,
 $\int \|x\|^p \mu  (dx) < \infty$. Then $W_p(\mu _n ,\mu ) \to 0$ as $n \to \infty $ is equivalent
to each of the following.

(a) $\mu _n\Longrightarrow \mu $ and $\int \|x\|^p \mu _n (dx) \to \int \|x\|^p \mu (dx)$.

(b) $\mu _n\Longrightarrow \mu $ and $\|\cdot\|^p$ is uniformly  $\mu _n$-integrable.

(c) $\int \phi (x) \mu _n (dx) \to \int \phi (x) \mu (dx)$ for every continuous $\phi $ such that $\phi (x)=O(\|x\|^p)$ at infinity.

\smallskip
Since in general an analog of the representation (2) does not exist for probability measures on a Banach 
space, Bickel and Freedman proved, in their setting, the existence of a solution to the transportation problem
for $c(x,y)=\|x-y\|^p$.

Recently, Ambrosio, Gigli, and Savar\'e 
 proved (\cite{Amb}, Proposition 7.1.5) an analog of part (b) of the above result for 
probability measures on a Radon space $X$ (see also Lemma 5.1.7 and
Remark 7.1.11 there).
These authors also established the existence of a solution 
to the Kantorovich transportation problem in $X$ for a wide class of cost functions.
We use this existence result
to prove criteria for the convergence in $T_c$ with $c(x,y)=C(d(x,y))$,
where $C$ is a non-decreasing continuous function satisfying the doubling condition (6) which controls 
the rate of growth of $C$ (Theorem 2 and Corollary 1). 
Since the class of such cost functions includes all the  $d^p$s, $p \geq 1$,  
the convergence results of Bickel and Freedman as well as those of  Ambrosio, Gigli, and Savar\'e follow from Corollary 1.
Note that instead of the Radon space $X$ (a separable metric space, where, by definition, every probability
measure is tight), we consider more familiar in the theory of probability object, a complete separable space $(M,d)$
where completeness and separability together provide the tightness of a probability measure; all our arguments
remain true for a Radon space (see also Remark 1).

In Theorem 2 we also obtain sufficient conditions for the convergence of probability measures 
in the transportation distance without assuming the doubling condition on $C$.
For instance, any convex $C:{\bf R}^+ \to {\bf R}^+$ with $C(0)=0$ satisfies Theorem 2.
We then provide an example
of a $C$ growing exponentially fast for which the converse
implication does not hold. 

\section{Convergence in Transportation Distance}

The following result of Ambrosio, Gigli, and Savar\'e  \cite{Amb} asserts the existence in  $\Pi$ of a probability measure  
which minimizes the total
transportation cost under rather weak assumptions on the
cost function. For the sake of completeness, we include a self-contained proof in Section 4.

\medskip
\noindent
{\bf Theorem 1.} {\it Let $(M,d)$ be a complete separable metric space, 
and let $T_{c}(\mu ,\nu )$ be defined by $(1)$
with $c:M \times M \to [0,+\infty)$ lower semicontinuous. 

Then there exists $\pi ^* \in \Pi$ such that $\int c(x,y)d\pi ^*(x,y)=T_{c}(\mu ,\nu )$.
Or, equivalently, there exists a pair of random elements $X$ and $Y$ with respective distributions $\mu $ and
$\nu$, such that ${\bf E}c(X,Y)=T_{c}(\mu ,\nu )$.}

\medskip
\noindent
{\bf Remark 1.} In the corresponding statement in \cite{Amb} the space $(M,d)$ need not be a complete
separable metric space but just a Radon space. In fact, our proof also shows that  completeness is unnecessary and 
that the tightness of $\mu $ and $\nu$ will suffice. On the other hand, the hypothesis of separability of $(M,d)$
can be weakened to the topological separability if both $\mu $ and $\nu $ have  separable supports (see Billingsley \cite{Bil2},
Appendix III).

\medskip
The Kantorovich problem is closely related to the  Monge transportation problem which is 
the problem of finding a map $s^*$
pushing $\mu $ forward to $\nu $ (i.e. such that $\nu (B)=\mu (s^{-1}(B))$ for any Borel set $B$) and minimizing the total transportation cost:
$\inf _s \int c(x,s(x))d\mu=\int c(x,s^*(x))d\mu,$ 
where the infimum is taken over all Borel maps $s$
pushing $\mu $ forward  to $\nu$.  A solution $s^*$ to the Monge transportation problem uniquely 
determines a probability measure $\pi ^*$ on $M \times M$ such that
the random elements $X$ and $Y$, $Y=s^*(X)$, with respective distributions $\mu $ and $\nu$ have  joint law $\pi ^*$.
This measure $\pi ^*$ minimizes the Monge transportation cost:
\begin{equation}
\int c(x,y)d\pi^* (x,y)=\inf _{\Pi^*} \int c(x,y)d\pi (x,y),
\end{equation}
where the infimum is taken over the set $\Pi^*$ of joint distributions of $M$-valued random elements
$X$ and $Y$ with respective distributions $\mu $ and $\nu $ and such that $Y$ is measurable with
respect to the Borel field $\sigma (X)$. Comparing (4) and (1) yields the relation
$\Pi^* \subset \Pi$, which  immediately leads
to the following conclusions: (i) the (Kantorovich) transportation distance $T_c(\mu, \nu)$
never exceeds the Monge transportation distance $\tilde T_c(\mu,\nu)$,
$$
\tilde T_c(\mu,\nu)=\inf _{\Pi^*} \int c(x,y)d\pi (x,y)=\inf _s \int c(x,s(x))d\mu;
$$
(ii) a probability measure $\pi ^*$ corresponding to the solution $s^*$ of the Monge transportation
problem (MTP) is not necessarily  a solution to the Kantorovich transportation problem (KTP);
conversely, a solution $\pi '$ of the KTP, where $\pi '$ is the joint distribution of $X$ and $Y$,
is a solution to the MTP if and only if there exists
a Borel map $s'$ such that $Y=s'(X)$. 

For random elements $X$ and $Y$ in a Hilbert space 
Cuesta and Matran \cite{Cuesta} have provided conditions for the existence of an increasing map $s$, $s(X)=Y$, 
such that $W_2^2(\mu, \nu)={\bf E}\|X- s(X)\|^2$, i.e. $X$ and $Y=s(X)$ give the solution
to both the MTP and the KTP. They also showed that if $\mu $ is either absolutely continuous 
with respect to the Lebesgue measure on ${\bf R}^k$ or is a Gaussian measure on a Hilbert space, then these conditions 
are satisfied. For compactly  supported absolutely continuous distributions 
on ${\bf R}^k$ and a convex cost function $c(x-y)$ 
Caffarelli \cite{Caf} has determined the form of the optimal map (the solution to the MTP)
as a gradient of $c$; the uniqueness of the solution is also obtained there.
Simultaneously, Gangbo and McCann \cite{GM} proved the same results
for non-necessarily boundedly supported probability measures. They
also showed that the solution to the MTP is the KTP solution as well, and that
a similar result holds true for $c(x,y)=l(\|x-y\|)$, where $l$ is  strictly  concave.
Note that in  all the  existence statements mentioned above,
the conditions of Theorem 1 are satisfied. A comprehensive review of the results on
the solutions to the KTP and the MTP  can be found in the books of Rachev and R\"uschendorf
\cite{R&R}.

The main result of the work presented here is now given.

\medskip
\noindent
{\bf Theorem 2.} {\it Let $\mu _n $ and $\mu $ be probability
measures on a complete separable metric space $(M,d)$ and let 
$c: M \times M \rightarrow {\bf R}$ be such that $c(x,y)=C(d(x,y))$,
where $C:[0,+\infty) \to [0,+\infty)$ is a non-decreasing continuous function with $C(0)=0$. Let
\begin{equation}
\int C(2d(x,a))\mu _n (dx) < \infty,\ \ \ \ \int C(2d(x,a)) \mu  (dx) < \infty
\end{equation}
for some (and, therefore, for all) $a \in M$.
Then} 
\begin{equation*}
\begin{rcase}
   (a)&\  \mu _n\Longrightarrow \mu,\\
   (b)& \int C(2d(x,a))\mu _n (dx) \rightarrow \int C(2d(x,a)) \mu  (dx)
\end{rcase}
\Longrightarrow T_c(\mu _n, \mu )\rightarrow 0.
\end{equation*}
{\it Conversely, if  $T_c(\mu _n, \mu )\rightarrow 0$, then}  $\mu _n\Longrightarrow \mu $.
{\it If, additionally, $C$ satisfies a doubling condition, i.e. if there exists
a positive constant  $\lambda $ such that
for all $y \geq 0$ }
\begin{equation}
C (2y) \leq \lambda C(y),
\end{equation}
{\it then}
\begin{equation*}
T_c(\mu _n, \mu )\rightarrow 0 \Longleftrightarrow
 \begin{cases}
   (a)\  \mu _n\Longrightarrow \mu,\\
   (b) \int C(2d(x,a))\mu _n (dx) \rightarrow \int C(2d(x,a)) \mu  (dx).
   \end{cases}
\end{equation*}
\medskip
\noindent
{\bf Corollary 1.} {\it If, in the setting of Theorem $2$, 
$C$ satisfies $(6)$, then 
$$
(6) \Longleftrightarrow  \int C(d(x,a))\mu _n (dx) < \infty, \int C(d(x,a)) \mu  (dx) < \infty ,
$$
and thus}
\begin{equation*}
T_c(\mu _n, \mu )\rightarrow 0 \Longleftrightarrow
 \begin{cases}
   (a)\  \mu _n\Longrightarrow \mu,\\
   (b') \int C(d(x,a))\mu _n (dx) \rightarrow \int C(d(x,a)) \mu  (dx).
   \end{cases}
\end{equation*}

\medskip
Corollary 1 is equivalent to a result of Rachev (Theorem 1 in \cite{Rachev2}) proved 
by using the relations between the L\'evy-Prokhorov metric and the $T_c$-distance.
Since for any $p\geq 1$, the function $c(x,y)=d^p(x,y)$ satisfies
the conditions of Theorem 2 as well as (6) with $\lambda =2^{p}$, 
Corollary 1 recovers part (a) in the result of Bickel and Freedman mentioned above. 
Ambrosio, Gigli, and Savar\'e  \cite{Amb} proved an analog of Theorem 2 in a Hilbert space when cost function is
continuous, strictly increasing and surjective map.

Note that the class of functions $C$ covered by Theorem 2 
includes functions with a faster than polynomial rate of growth at infinity (e.g.
 $C(d(x,y))=\exp(d(x,y))-1$).
For functions $C$ growing exponentially fast at infinity, and in contrast to $C(d(x,y))=d^p(x,y)$,
$T_c(\mu _n, \mu )\rightarrow 0$  need not imply the convergence of 
$\int C(2d(x,a))\mu _n (dx)$ to $\int C(2d(x,a)) \mu  (dx)$, for some $a \in M$. Indeed, 
one can take
the probability measures $\mu _n$ and $\mu $ on ${\bf R}$ defined in Example 1, below, and 
$c(x,y)=C(|x-y|)=\exp (|x-y|)-1$. 

\medskip
As a corollary to Theorem 2 we obtain the following result relating
the convergence in total variation to the
convergence in transportation distance. It is well known 
that the total variation distance itself is a particular case of  
transportation distance (with $c(x,y)=2{\bf 1}_{\{x\neq y\}}$). 

\medskip
\noindent
{\bf Corollary 2.} {\it Let $\mu $ and $\nu $ be boundedly supported probability measures on
a complete separable metric space $(M,d)$, and let
$\phi $ be a continuous function
on  $(M,d)$. Then
\begin{equation*}
\left |\int \phi(x)\mu (dx)-\int \phi (x)\nu (dx) \right |\leq L _{\phi} \|\mu -\nu \|_{TV}
\end{equation*}
for some positive constant  $L _{\phi}$.

Let $\mu _n$ be  probability measures on $M$ with respective supports $K_n$, $n \geq 1$. Let $\cup _nK_n$ be bounded.
If $c(x,y)=C(d(x,y))$, where 
$C:[0,+\infty) \to [0,+\infty)$ is  non-decreasing, continuous with $C(0)=0$,
then} 
$$ \|\mu_n -\mu \|_{TV} \rightarrow 0 \Rightarrow T_c(\mu_n, \mu )\rightarrow 0.$$

\medskip
Without the boundedness restriction on
$\cup K_n$ the last implication is not true, as the following example shows.  

\medskip
\noindent
{\bf Example 1.} Let $\mu $ be the uniform distribution
on $(0,1)$ and, for all $n\in {\bf N}$, let 
$$\mu _n (dx)=\frac {n-1}{n}{\bf 1}_{(0,1)}(x)dx+ \frac 1n \delta _{x_n}(dx),
$$
$x_n \notin (0,1)$. Then
$$\|\mu_n -\mu \|_{TV}=\int _0^1 |f_n(x)-1|dx+\mu _n(x_n)=
\frac 2n\rightarrow 0$$
as $n\rightarrow \infty $. Hence $\mu _n
\stackrel {TV}{\longrightarrow } \mu$
for any choice of the sequence $(x_n)$.
Let $c(x,y)=C(|x-y|)$, with  $C:[0,+\infty) \to [0,+\infty)$, $C(0)=0$,  convex,
also satisfying (6). Then,
$$\int C(|x|)\mu (dx)=\int _0^1 C(|x|)dx\leq \max _{0< |x| < 1}C(|x|)<+\infty,$$
$$
\int C(|x|)\mu _n (dx)=\int _0^1\frac {n-1}{n}C(|x|)dx+ \mu _n(x_n)
C(|x_n|)\leq \max _{0< |x| < 1}C(|x|)\frac {n-1}{n} +\frac {C(|x_n|)}{n}<+\infty, $$
for any $n$. So  all the conditions of Corollary 1
are satisfied. 
Since weak convergence is implied by convergence in total variation,
$T_c(\mu _n,\mu)\rightarrow 0$ holds
if and only if $\int C(|x|)\mu _n (dx)\rightarrow \int C(|x|)\mu (dx)$.
Take $x_n=2^n$, then $C(|x_n|)=C(2^n)\geq 2^{n-1}C(2)$ and
$C(|x_n|)/n \rightarrow +\infty$ as $n \to \infty$. 
Therefore, 
$$ \int C(|x|)\mu _n (dx)\geq \frac {C(|x_n|)}{n} \rightarrow +\infty \neq \int C(|x|)\mu (dx) .$$
By Corollary 1, $T_c(\mu _n,\mu )$ does not converge to $0$.

\section{Applications to the Central Limit Theorem}

Next, we apply Theorem 2 to
obtain the CLT in the transportation distance.
We provide  sufficient conditions for the convergence
of the laws of the normalized sums to the standard Gaussian measure on ${\bf R}$
for stationary sequences which are either independent, strongly mixing or associated.

\subsection{Independent sequences}

Let $(X_n)$ be a sequence of
independent identically distributed random variables, ${\bf E}X_1=0$,
${\bf E}X_1^2=\sigma ^2$, $0<\sigma <+\infty$. Let $S_n=\sum _{k=1}^n X_k$. 
Then by the classical central limit theorem
$$\frac {S_n}{\sigma \sqrt n}\stackrel d {\longrightarrow }
Z \sim N(0,1).$$

Let $\mu _n$ denote the probability law of the normalized sum
$S_n/(\sigma \sqrt n)$, and let $\gamma $ be the standard Gaussian measure on ${\bf R}$. 
We find additional conditions on the sequence $(X_n)$ and
on the cost function
to obtain the convergence of $\mu _n$ to $\gamma $ 
in the $T_c$-distance.

\smallskip
\noindent
{\bf Theorem 3.} {\it Let $c(x,y)=C(|x-y|)$,
where $C:[0,+\infty) \to [0,+\infty)$ is a non-decreasing continuous 
function with} $C(0)=0$.

(i) {\it If there exists $p \geq 2$ such that $C(x)=O(x^p) $ at infinity and 
${\bf E}|X_1|^p<+ \infty $,
then} $T_c(\mu_n, \gamma ) \to 0$.

(ii) {\it Otherwise, let ${\bf E}C(4 \sqrt 2|Z|) <+\infty$ and
let $\sum _{k=1}^{\infty }k^k {\bf E}X_1^{2k}<+\infty$,  then}
$T_c(\mu_n, \gamma ) \to 0$.

\smallskip
The CLT in the $W_2$-distance was proved by Tanaka \cite{Tanaka} for distributions 
on ${\bf R}$ and by Cuesta and Matran \cite{Cuesta} for distributions on a Hilbert space;
both results require the finiteness of the fourth moment.
Very recently, Johnson and Samworth \cite{JSam}, \cite{JSam1} proved that $W_p(\mu _n,\gamma )\rightarrow 0$,
$p \geq 2$, under the condition ${\bf E}|X_1|^p< \infty $. This statement is a particular case of 
part (i) of Theorem 3. However, these authors also proved the convergence to an $\alpha $-stable law
in the Wasserstein-Mallows $\alpha$-distance, $\alpha  \in (0,2)$. 

We will prove Theorem 3 by applying Theorem 2.
The CLT yields  weak convergence; therefore, to prove convergence in 
$T_c$-distance, we need to verify the convergence of 
$\int C(2|x|)d\mu _n $ to 
$\int C(2|x|)d\gamma $. To do so, we use
Rosenthal's inequality  which asserts that
for stationary independent sequence $(X_n)$ of centered random variables
\begin{equation}
{\bf E}|S_n|^p \leq K(p) \left(n{\bf E}|X_1|^p+ n^{p/2} ({\bf E}|X_1|^2)^{p/2} \right)
\end{equation}
for $p>1$ and a positive constant $K(p)$ depending only on $p$ (Petrov \cite{Petrov}).

\subsection{Strong mixing sequences}

The coefficients $\alpha _n $ 
of strong mixing of a random sequence $(X_n)$ are defined as
$$
\alpha _n=\sup\{|P(A\cap B)-P(A)P(B)|:A\in {\mathcal F}_1^k,\ 
B \in {\mathcal F}_{k+n}^{+\infty }, k\geq 1\},
$$
where ${\mathcal F}_k^{k+m}$ is the $\sigma $-field
generated by the random variables $X_k,X_{k+1},...,X_{k+m}$.
A sequence is said to satisfy a strong mixing
condition if $\alpha _n \rightarrow 0$
as $n\rightarrow +\infty $. 

A CLT for a stationary strong mixing sequence $(X_n)$ was proved by Denker
\cite{Den} in the following form.
Let  ${\bf E}X_1=0$,
${\bf E}X_1^2=\sigma ^2$, $0<\sigma <+\infty$, and $\sigma _n^2={\bf E}S_n^2=
nh(n)$, where $h(n)$ is a slowly varying function.
Let $S_n=\sum _{k=1}^n X_k$, and let $\mu _n$ be the law of $S_n/\sigma_n$.
Then $\mu _n \Longrightarrow \gamma$, where $\gamma $ is the standard Gaussian measure on ${\bf R}$.

To obtain the convergence of $\mu _n$ to $\gamma $ in 
$T_c$, we need  additional conditions on the rate of
decay of the coefficient $\alpha _n$ providing a moment
inequality for sums. Such a result exists, it is due to 
Yokoyama \cite{Y} and asserts that if $(X_n)$ is a 
stationary strong mixing sequence such that
${\bf E}X_1=0$, ${\bf E}|X_1|^{p+\delta }<+\infty $, $p>2$, $\delta >0$
and
\begin {equation}
\sum _{n=1}^{\infty }(n+1)^{\frac p2 -1}(\alpha _n)^{\frac {\delta}
{p+\delta }}<+\infty,
\end {equation}
then 
\begin{equation}
{\bf E}|S_n|^p\leq K(p)n^{\frac p2},
\end{equation}
where the positive constant $K(p)$ depends only on $p$.

\medskip
\noindent
{\bf Theorem 4.} {\it Let $c(x,y)=C(|x-y|)$,
where $C:[0,+\infty) \to [0,+\infty)$ is a non-decreasing continuous 
function with} $C(0)=0$.

(i) {\it If there exist $p > 2$ and $ \delta >0$ such that the condition $(8)$ is satisfied
and $C(x)=O(x^p) $ at infinity,
then} $T_c(\mu_n, \gamma ) \to 0$.

(ii) {\it Otherwise, let ${\bf E}C(4 \sqrt 2|Z|) <+\infty$,
let $\sum _{k=1}^{\infty }k^k {\bf E}X_1^{2k}<+\infty$, and let $(8)$ be satisfied for all $p>2$, then}
$T_c(\mu_n, \gamma ) \to 0$.

\smallskip
To prove Theorem 4, we once again apply Theorem 2.
Since the corresponding CLT implies weak convergence to the Gaussian measure,
it is sufficient to show the convergence of  
 ${\bf E}C(2|Y_n|)$ to ${\bf E}C(2|Z|)$.
 
\subsection{Associated sequence}

Recall that a set of random variables $\xi
=\left( \xi _1,...,\xi _m\right) $ is called associated if
for any two coordinatewise increasing functions $f,\ g:{\bf R}^m\rightarrow {\bf R}$ 
such that ${\bf E}f(\xi )g(\xi )$, ${\bf E}f(\xi )$ and ${\bf E}g(\xi )$ 
exist, 
$$\mbox{Cov} \left( f(\xi ),g(\xi )\right) \geq 0. $$
An infinite set of random variables is  associated if
all of its finite subsets are associated.

Newman \cite{N80} proved the CLT for associated sequence under the following conditions. Let $(X_n)$ be a 
stationary associated sequence,
${\bf E}X_1=0$,
${\bf E}X_1^2=\sigma ^2$, $0<\sigma <+\infty$, $\sigma _n^2={\bf E}S_n^2=
nh(n)$, where $h(n)$ is a slowly varying function.
Let $\mu_n$ be the law of $S_n /\sigma _n$ and let $\gamma $ be the standard Gaussian measure on ${\bf R}$. 
Then $\mu_n \Longrightarrow \gamma$.

Asymptotic independence
for associated sequence $(X_n)$ is usually stated in
terms of the Cox-Grimmett coefficient $u(n)$ defined by:
$$ u(n)=\sup_{k \geq 1}\sum_{j: |j-k| \geq n} \mbox{Cov}(X_j,X_k). $$
For a stationary sequence the Cox-Grimmett coefficient is just
the tail of the series of covariances:
$$ u(n)=2\sum _{k=n+1}^{\infty }\mbox{Cov}(X_1,X_k).$$

To prove the convergence of $\mu _n$ to $\gamma $ in the transportation
distance, we use a condition on the rate of
decay of the Cox-Grimmett coefficient. This condition implies the following moment
inequality for sums (Birkel \cite{B}). If $(X_n)$ is a 
stationary associated sequence,
${\bf E}X_1=0$, ${\bf E}|X_1|^{p+\delta }<+\infty $, $p>2$, $\delta >0$
and
\begin {equation}
u(n)\leq Bn^{-\frac {(p-2)(p+\delta )}{2\delta }},
\end {equation}
then 
\begin{equation}
{\bf E}|S_n|^p\leq K(p)n^{\frac p2},
\end{equation}
where the positive constant $K(p)$ depends only on $p$.

\medskip
\noindent
{\bf Theorem 5.} {\it Let $c(x,y)=C(|x-y|)$,
where $C:[0,+\infty) \to [0,+\infty)$ is a non-decreasing continuous 
function with} $C(0)=0$.

(i) {\it If there exist $p > 2$ and $ \delta >0$ such that the condition $(10)$ is satisfied
and $C(x)=O(x^p) $ at infinity,
then} $T_c(\mu_n, \gamma ) \to 0$.

(ii) {\it Otherwise, let ${\bf E}C(4 \sqrt 2|Z|) <+\infty$,
let $\sum _{k=1}^{\infty }k^k {\bf E}X_1^{2k}<+\infty$, and let $(10)$ be satisfied for all $p>2$, then}
$T_c(\mu_n, \gamma ) \to 0$.

\section{Proofs}
\begin{proof}[{\bf Proof of Theorem 1.}] 
We first show that $\Pi $ is a tight set. Indeed, for any positive $\varepsilon $ there exist
compact sets $K_1$, $K_2 \in \mathcal B (M)$, such that $\mu (K_1) \geq 1- \frac {\varepsilon }{2}$
and $\nu (K_2) \geq 1- \frac {\varepsilon }{2}$.
Let $\pi \in \Pi$ and let $(X,Y)$ be a random vector  with law $\pi $. Then,
\begin{equation}
   \begin{aligned}
   \pi (K_1 \times K_2)&=P(X\in K_1,Y\in K_2)=P(X\in K_1)+P(Y\in K_2)-P((X\in K_1)\cup (Y\in K_2))\\
                       &\geq  \mu (K_1)+ \nu (K_2) -1 \geq (1- \varepsilon /2)+(1- \varepsilon /2)-1=1- \varepsilon.
		    \end{aligned}
\end{equation}
Since (12) holds for all $\pi \in \Pi$ with the same compact set $K_1 \times K_2$, this proves that $\Pi $ is tight.
Therefore, according to Prokhorov's theorem (Billingsley \cite{Bil2}, Section 5), $\Pi $ is relatively compact. 

If $T_{c}(\mu ,\nu )=+\infty$, then $\int c(x,y)d\pi (x,y)=+\infty $, for all $\pi \in \Pi$ and
$\pi ^* $ can be chosen to be any probability measure from $\Pi$. 

If $T_{c}(\mu ,\nu )< +\infty$, then there exists a sequence $\pi _n$ from $\Pi $ such that
\begin{equation}
\int c(x,y)d\pi _n(x,y) \to T_{c}(\mu ,\nu ).
\end{equation}
On the other hand, the relative compactness of $\Pi $ implies the existence of
a subsequence $\pi _{n_k}$
which converges weakly to some probability measure $\pi $ on $\mathcal B(M \times M)$.
Let us verify that $\pi $ is the measure $\pi ^*$ we are looking for.
First, we want to prove that $\pi \in \Pi$, i.e. that the marginal distributions of $\pi $
are $\mu$ and $\nu$, respectively. 

Let $\mu_1$ and $\nu _1$ be marginals of $\pi $. We will check that $\mu_1(B)=\mu (B)$,
for any $B \in \mathcal B (M)$ such that $\mu _1(\partial B)=0$. Indeed,
since $\partial (B \times M) \subset (\partial B \times M) \cup (B \times \partial M)
=\partial B \times M$
(Billingsley \cite{Bil2}, (2.8)), we have
$$
\pi (\partial (B \times M)) \leq \pi (\partial B \times M)=\mu _1(\partial B)=0.
$$
Therefore, the weak convergence  $\pi _{n_k} \Rightarrow \pi $ implies that $\pi _{n_k}(B \times M) 
\to \pi (B \times M)$,
and we obtain
\begin{equation*}
  \mu (B)=\pi _{n_k}(B \times M) \to \pi (B \times M)=\mu _1(B).
\end{equation*}
Similarly, we can show that $\nu_1(B)=\nu (B)$,
for any $B \in \mathcal B (M)$ such that $\nu _1(\partial B)=0$. Finally, it remains to check that
two probability measures $\mu _1$ and $\mu $ (respectively $\nu _1$ and  $\nu$) are the same if they coincide
on the Borel sets having a boundary of $\mu _1$-measure (respectively $\nu_1$-measure) zero. 

Let $D\in \mathcal B (M)$ be a closed set. For $\varepsilon >0$, let
$D^{\varepsilon }=\{x \in M: d(x,D) < \varepsilon \}$
and let ${\mathcal D}= \{D^{\varepsilon }, 0 <\varepsilon <1\}$. Then  there
exists at most a countable number of $\varepsilon _k$,
$0<\varepsilon _k<1$, such that sets $D^{\varepsilon _k}$ have
a boundary of positive $\mu _1$-measure. We remove the sets $D^{\varepsilon _k}$ from ${\mathcal D}$,
and obtain  
$$
{\mathcal D }^0=\{D^{\varepsilon }, 0 <\varepsilon <1, \mu_1(\partial D^{\varepsilon })=0\}.
$$
We can then choose a decreasing sequence $\varepsilon _n \to 0$, $0<\varepsilon _n<1$, 
with $D_n=D^{\varepsilon _n}\in {\mathcal D}^0$.
The sets $D_n$ are such that: (a) $D_{n+1} \subset D_{n}$ for all $n$; (b) $\bigcap _n D_n = D\cup \partial D=D$;
(c) $\mu_1(D_n)=\mu (D_n)$.
The properties (a)--(c) yield
$$
\mu _1(D)= \mu _1(\bigcap _n D_n)= \lim _{n \to \infty} \mu_1 (D_n)=\lim _{n \to \infty} \mu (D_n)
=\mu (D).
$$
Therefore, the measures $\mu_1 $ and $\mu $ coincide on all the closed subsets of $M$.
Since $\mathcal B(M)$ is generated by such sets, we conclude that $\mu_1=\mu$. Similar arguments lead to $\nu_1=\nu$. We have proved that
the probability measure $\pi$ has respective marginals $\mu$ and $\nu$, i.e. $\pi \in \Pi$. 

Next, we will check that $\int c(x,y)d\pi (x,y)=T_{c}(\mu ,\nu )$. Since $c$ is lower semicontinuous, 
for any real $b$ the set $\{(x,y): c(x,y)>b\}$ is open (\cite{Bil2}, Appendix I). Let 
$A=\{(x,y): c(x,y)>0\}$. Then the weak convergence $\pi _{n_k} \Longrightarrow \pi$ and (13) 
imply that
\begin{equation*}
\begin{aligned}
\int c(x,y)d\pi (x,y)&=\int _{A}c(x,y)d\pi (x,y)\\
                     &\leq \lim\ \inf_{n_k} \int _{A}c(x,y)d\pi _{n_k}(x,y)\\
                     &=\lim\ \inf_{n_k} \int c(x,y)d\pi _{n_k}(x,y)\\
		     &=T_{c}(\mu ,\nu ).
		     \end{aligned}
		     \end{equation*}
Since $\pi \in \Pi$, the reverse inequality $\int c(x,y)d\pi (x,y) \geq T_{c}(\mu ,\nu )$ holds true.
We thus conclude that $\int c(x,y)d\pi (x,y) = T_{c}(\mu ,\nu )$. In other words, the transportation distance
becomes the total transportation cost associated to the measure $\pi $.
Finally, we set $\pi ^*=\pi $ and
the proof is now complete.
\end{proof}

\smallskip
\noindent
\begin{proof}[{\bf Proof of Theorem 2 and Corollary 1.}] Assume that both (a) and (b) are satisfied. 
Let $X$, $X_n$ be
random elements with respective distributions $\mu $
and $\mu _n$ and such that $X$ and $X_n$ are independent, for any $n$. Then $(C(2d(X_n,a)))$ is uniformly bounded,  that is
$I_1= \sup _n {\bf E}C(2d(X_n,a))<\infty.$
Set $I_2={\bf E}C(2d(X,a))<\infty$.

Fix $\varepsilon >0$ and choose a compact set $K_1$ in ${\mathcal B}(M)$ such that
$\mu (\partial K_1)=0$ and
$$
\int _{(K_1)^c}C(2d(x,a))d\mu (x) <\varepsilon.
$$
The weak convergence $\mu_n \Longrightarrow \mu$ implies the tightness of
the family $(\mu_n, \mu )_{n \geq 1}$, thus there exists a compact set $K_2 \in {\mathcal B}(M)$ such that
$\mu_n(K_2)^c< \varepsilon$, $\mu (K_2)^c< \varepsilon$ and $\mu (\partial K_2)=0$.
Let $K=K_1 \cup K_2$. Then  $K$ is  compact, and
\begin{equation}
\int _{K^c}C(2d(x,a))d\mu (x) <\varepsilon,
\end{equation}
\begin{equation}
  \begin{array}{ll}
  \mu_n(K^c)< \varepsilon,&\mu (K^c)< \varepsilon,
  \end{array}
\end{equation}
with also $\mu (\partial K)=0$, since $\mu (\partial K)\leq \mu (\partial K_1)+\mu (\partial K_2)$.
Since (b) holds, we can choose a positive integer $N_1$ such that for any $n\geq N_1$,
\begin{equation}
\left|\int C(2d(x,a))d\mu _n (x)- \int C(2d(x,a))d\mu (x)\right|<\varepsilon.
\end{equation}
 As $X_n\stackrel d {\longrightarrow} X$,
 for the chosen compact set $K$ and the continuous function $C(2d(\cdot,a))$ we have
$${\bf E}C(2d(X_n,a)){\bf 1}_{\{X_n\in K\}}\rightarrow 
{\bf E}C(2d(X,a)){\bf 1}_{\{X\in K\}}.$$
Hence, we can choose a positive integer $N_2$ such that, for any $n\geq N_2$,
\begin{equation}
\left|\int _K C(2d(x,a))d\mu _n (x)- \int _K C(2d(x,a))d\mu (x)\right|<\varepsilon.
\end{equation}
Then for $n\geq \max \{N_1,N_2\}$, the estimates (14), (16) and (17) yield
\begin{equation}
\begin{aligned}
\left|\int _{K^c} \right.& \left. C(2d(x,a))d\mu _n (x)\right|\\
                            &\leq \left|\int C(2d(x,a))d\mu _n (x)- \int C(2d(x,a))d\mu (x)\right|\\
			    &+\left|\int _K C(2d(x,a))d\mu _n (x)- \int _K C(2d(x,a))d\mu (x)\right|+
\left|\int _{K^c}C(2d(x,a))d\mu (x)\right|\\
    &<3 \varepsilon.
\end{aligned}
\end{equation}
The weak convergence  $X_n\stackrel d {\longrightarrow} X$ implies that
$C(2d(X_n,X)){\bf 1}_{\{X_n\in K,\ X\in K\}}\stackrel d {\longrightarrow} 0.$
The continuous function $C(d(x,y))$ is bounded on the compact set $K \times K$, therefore
$${\bf E}C(d(X_n,X)){\bf 1}_{\{X_n\in K,\ X\in K\}}
\rightarrow 0.$$
This means that there exists a positive integer $N_3$ such that, for any $n\geq N_3$, 
\begin{equation}
\left | \int _K\int _K C(d(X_n,X))d\pi _n (x,y)\right|<\varepsilon,
\end{equation}
where $\pi _n$ is the joint distribution of $X_n$ and $X$.

Since $C$ is a non-negative and non-decreasing, 
\begin{equation}
\begin{aligned}
C(d(x,y)) &\leq C(d(x,a)+d(y,a))\leq  C(2 \max \{d(x,a),d(y,a)\})\\
          &\leq C(2d(x,a))+C(2d(y,a)),
	  \end{aligned}
\end{equation}
for all $x,y \in M$.

Using the inequalities (14), (15), (18), (20), and the independence of $X_n$ and $X$, we have:
\begin{equation}
   \begin{aligned}
    \int _{K^c}&\int _{K^c} C(d(x,y))d\pi _n (x,y)={\bf E}C(d(X_n,X)){\bf 1}_{\{X_n \in K^c,\ X \in K^c\}}\\
    &\leq {\bf E}C(2d(X_n,a)){\bf 1}_{\{X_n \in K^c\}} {\bf 1}_{\{X \in K^c\}}+ {\bf E}C(2d(X,a)){\bf 1}_{\{X \in K^c\}}{\bf 1}_{\{X _n\in K^c\}}\\
    &\leq \left(\int _{K^c} C(2d(x,a))d\mu _n (x) \right)\mu (K^c) + \left(\int _{K^c} C(2d(x,a))d\mu (x) \right)\mu _n (K^c)\\
    &< 3\varepsilon ^2+  \varepsilon^2,
    \end{aligned}
    \end{equation}
for all $n \geq \max \{N_1, N_2\}$. Similarly, 
\begin{equation}
   \begin{aligned}
    \int _{K}&\int _{K^c} C(d(x,y))d\pi _n (x,y)={\bf E}C(d(X_n,X)){\bf 1}_{\{X_n \in K,\ X \in K^c\}}\\
    &\leq {\bf E}C(2d(X_n,a)){\bf 1}_{\{X_n \in K\}} {\bf 1}_{\{X \in K^c\}}+ {\bf E}C(2d(X,a)){\bf 1}_{\{X \in K^c\}}{\bf 1}_{\{X _n\in K\}}\\
    &\leq I_1 \mu (K^c) + \varepsilon \mu _n (K)\\
    &< I_1 \varepsilon +  \varepsilon,
    \end{aligned}
    \end{equation}
    and
\begin{equation}
   \begin{aligned}
    \int _{K^c}&\int _{K} C(d(x,y))d\pi _n (x,y)={\bf E}C(d(X_n,X)){\bf 1}_{\{X_n \in K^c,\ X \in K\}}\\
    &\leq {\bf E}C(2d(X_n,a)){\bf 1}_{\{X_n \in K^c\}} {\bf 1}_{\{X \in K\}}+ {\bf E}C(2d(X,a)){\bf 1}_{\{X \in K\}}{\bf 1}_{\{X _n\in K^c\}}\\
    &\leq 3\varepsilon \mu (K) + I_2 \mu _n (K^c)< 3 \varepsilon + I_2 \varepsilon,
    \end{aligned}
    \end{equation}
for $n \geq \max \{N_1, N_2\}$.

Thus for $n\geq \max \{N_1,N_2,N_3\}$ the inequalities (19), (21)--(23) yield
\begin{equation*}
 \begin {aligned}
  T_c(\mu _n,\mu )&\leq {\bf E}C(d(X_n,X))=\int \int C(d(x,y))d\pi _n (x,y)\\
                  &=\int _K \int _K C(d(x,y))d\pi _n (x,y)+\int _{K^c}\int _{K^c} C(d(x,y))d\pi _n (x,y)\\
		  &+\int _{K}\int _{K^c} C(d(x,y))d\pi _n (x,y)+\int _{K^c}\int _{K} C(d(x,y))d\pi _n (x,y)\\
		   &\leq \varepsilon (6 +4 \varepsilon +I_1+I_2).
  \end{aligned}
\end{equation*}
We conclude that (a) and (b) imply $T_c(\mu _n,\mu )\rightarrow 0.$

Next, we assume that $T_c(\mu _n,\mu )\rightarrow 0$ and verify 
that (a) $\mu _n \Longrightarrow \mu $ takes place.
According to Theorem 1, for any $n$ there exists a pair of
random elements $X_n$ and $X$ with distributions $\mu _n$ and $\mu $,
respectively, which are minimizers of the total transportation cost:
$T_c(\mu _n,\mu )={\bf E}C(d(X_n,X)).$
Let us note that $X$ may depend on $n$, so each time it appears 
in this proof, we assume that $X=X^{(n)}$. (Of course all the $X^{(n)}$ have
the same law $\mu $.)

Since $C$ is a non-negative function, 
${\bf E}C(d(X_n,X))\rightarrow 0$, that is
$C(d(X_n,X))\stackrel {L_1}{\longrightarrow }0$. This implies that
\begin{equation}
C(d(X_n,X))\stackrel P {\longrightarrow }0. 
\end {equation}

Fix $\varepsilon >0$. As $C$
is non-decreasing, we have
$$\{d(X_n,X)>\varepsilon \}\subset \{C(d(X_n,X))\geq C(\varepsilon)\}.$$
The convergence result (24) means that the probability of the last 
event tends to 0, for any positive $C(\varepsilon)$, as $n \to \infty$. Hence for any
$\varepsilon >0$,
$P(d(X_n,X)>\varepsilon )\rightarrow 0,$
as $n\rightarrow \infty$. 
The convergence, in probability, of $d(X_n,X)=d(X_n,X^{(n)})\stackrel P {\longrightarrow }0$ implies
that $\mu _n\Longrightarrow \mu $ (Billingsley \cite{Bil2}, theorem 4.1).

Now, assume that the doubling condition (6) is satisfied 
and let us verify that $T_c(\mu _n,\mu)\to 0$ implies (b$'$).
Since $\mu _n \Longrightarrow \mu $ and since $C(d( \cdot, a))$ is continuous on $M$,
weak convergence holds: 
$C(d(X_n,a)\stackrel d {\longrightarrow }C(d(X,a))$. In order to
verify (b$'$), it thus suffices to check that the sequence $(C(d(X_n,a)))$
is uniformly integrable. The uniform integrability is equivalent to the pair of conditions: (i) $({\bf E}C(d(X_{n},a)))$
is uniformly bounded and (ii) for $A \in \mathcal F$,  $({\bf E}C(d(X_n,a)){\bf 1}_A)$ is uniformly continuous,
(i.e. $\sup  _{n} {\bf E}C(d(X_n,a)){\bf 1}_A \to 0$ as $P(A) \to 0$).

Together (6) and (20)  yield the 
inequalities
\begin{equation}
C(d(x,a)) \leq \lambda C \left( \frac 12 d(x,a) \right) \leq  \lambda C(d(x,y)) + \lambda C(d(y,a)),
\end{equation}
for all $x,y \in M$ and the positive constant $\lambda $. Then
\begin{equation}
{\bf E}C(d(X_{n},X))\geq \frac 1 {\lambda }{\bf E}C(d(X_{n},a))-{\bf E}C(d(X,a)).
\end{equation}
Suppose that $({\bf E}C(d(X_{n},a)))$ is not uniformly bounded.
Then there exists a subsequence $({\bf E}C(d(X_{n'},a)))$ such that
${\bf E}C(d(X_{n'},a))\rightarrow +\infty. $ Applying (26) to this
subsequence, we come to the following contradiction:
$$
0\leftarrow {\bf E}C(d(X_{n},X))\geq \frac 1 {\lambda }{\bf E}C(d(X_{n'},a))-{\bf E}C(d(X,a))
\rightarrow +\infty.
$$
Thus, $({\bf E}C(d(X_{n},a)))$ is uniformly bounded.

Let $\varepsilon $ be fixed, and let $A \in \mathcal F$. 
Since $T_c(\mu _n,\mu)\to 0$, we can choose a positive integer $N$ such that
${\bf E}C(d(X_{n},X)){\bf 1}_A < \varepsilon $, for all $n \geq N$. 
By applying once again the inequality (26), we obtain
$$
\sup _{n\geq N}{\bf E}C(d(X_n,a)){\bf 1}_A\leq
 \lambda \sup _n{\bf E}C(d(X_{n},X)){\bf 1}_A+ \lambda {\bf E}C(d(X,a)){\bf 1}_A.
$$
Let $P(A)\rightarrow 0$. Since $(C(d(X_{n},X)))$ is uniformly integrable
and since ${\bf E}C(d(X,a)) \leq {\bf E}C(2d(X,a))< \infty$,
$$\sup _{n\geq N}{\bf E}C(d(X_n,a)){\bf 1}_A\rightarrow 0,$$ 
i.e. $({\bf E}C(d(X_n,a)){\bf 1}_A)$ is uniformly continuous. Hence, 
the sequence $(C(d(X_n,a)))$
is uniformly integrable and (b$'$) $\int C(d(x,a))d\mu _n\rightarrow \int C(d(x,a))d\mu $ holds.

Note that from (6) and since $C$ is non-decreasing, the following two
inequalities hold true
\begin{equation*}
C(2d(x,a)) \leq \lambda C(d(x,a)),\ \ \ \ C(d(x,a)) \leq C(2d(x,a)),
\end{equation*}
for any $x \in M$. This implies that\\ 
$ (\int C(d(x,a))\mu _n (dx) < \infty)\Longleftrightarrow (\int C(2d(x,a))\mu _n (dx) < \infty)$
and that $ (\int C(d(x,a))\mu  (dx) < \infty)\Longleftrightarrow (\int C(2d(x,a))\mu  (dx) < \infty)$.
Therefore, in the setting of the theorem, the sequences $(C(d(X_n,a)))$ 
and $(C(2d(X_n,a)))$ are both either uniformly integrable or not,
and (b)$\Longleftrightarrow $ (b$'$).

This observation completes the proof of Theorem 2 and of Corollary 1.
\end{proof}

\smallskip
\noindent
\begin{proof}[{\bf Proof of Corollary 2.}] Let $K_1$ and $K_2$ be the respective supports
of $\mu $ and $\nu $. 
If $\mu $ and $\nu $ are absolutely continuous with respective densities
$f_1$ and $f_2$, then
\begin{equation*}
     \begin{aligned}
    \left |\int \phi (x)d\mu -\int \phi(x)d\nu \right |&= \left |\int \phi (x)f_1(x)d x -\int \phi (x)f_2(x)dx \right |\\
&\leq \int _{K_1\cup K_2} |\phi (x)| |f_1(x)-f_2(x)|dx\\
&\leq L_{\phi }\|\mu -\nu \|_{TV},
\end{aligned}
\end{equation*}
where $L_{\phi }=\sup \{|\phi (x)|:x\in \overline{(K_1\cup K_2)}\}$ (here $\overline A=A \cup \partial A$).

To prove the result in the general case, define the partition $(A_m)_{m \in {\bf Z}}$
of $M$, $A_m \in \mathcal B(M)$, as follows: 
$$
A_m=\{x \in M:m-1\leq \phi(x) < m\}.
$$
Thus,
\begin{equation}
\begin{aligned}
\left |\int \phi (x)d\mu -\int \phi (x)d\nu \right |&=
\left |\sum \limits_{m=-\infty}^{+\infty }\int \phi (x){\bf 1}_{A_m}d(\mu -\nu )\right |\\
 &\leq \sum \limits _{m=-\infty}^{+\infty }|m\|\mu (A_m)-\nu (A_m)|\\
 &\leq L_{\phi }\|\mu -\nu \|_{TV},
\end{aligned}
\end{equation}
where $L_{\phi }$ is defined as above, and where we used
the dual definition of the  total variation distance.
 
Let $\mu _n$ and $\mu $ be probability measures on
$M$ with  bounded supports respectively denoted $K_n$ and $K$. Let also
$\cup _n K_n$ be bounded and
$ \|\mu_n -\mu \|_{TV} \rightarrow 0$.
Convergence in total variation implies
weak convergence $\mu_n \Longrightarrow \mu $. All the conditions of Theorem 2 are satisfied.
Therefore, to prove the convergence of $\mu _n$ to $\mu $ in
$T_c$, it suffices to check  that 
 $\int C(2d(x,a))d\mu _n \rightarrow \int C(2d(x,a))d\mu $ for some $a \in M$.
The inequality (27) yields for any $n$
\begin{equation}
\left |\int \phi (x)d\mu _n -\int \phi (x)d\mu \right |\leq
L_{\phi }\|\mu _n-\mu\|_{TV},
\end{equation}
where $L_{\phi }=\sup \{|\phi (x)|:x\in \overline{\cup K_n}\}< \infty $ does not
depend on $n$. By fixing $a \in M$ and applying (28) to $\phi (x)=C(2d(x,a))$,  we obtain that the
convergence in total variation implies the convergence
of the integrals $\int \phi (x)d\mu _n \to \int \phi (x)d\mu$. This completes the proof.
\end{proof}

\smallskip
\noindent
\begin{proof}[{\bf Proof of Theorem 3.}] (i) First, we set $Y_n= S_n /(\sigma \sqrt n)$ and
check that ${\bf E}|Y_n|^p \to {\bf E}|Z|^p$. The classical CLT gives 
$Y_n \stackrel d \longrightarrow Z$, while the uniform boundedness of ${\bf E}|Y_n|^p$
follows from Rosenthal's inequality (7):
\begin{equation}
{\bf E}|Y_n|^p \leq K(p) \left(\frac  {{\bf E}|X_1|^p}{\sigma ^p n^{\frac p2-1}}+ 1 \right)
\leq K(p) \left(\frac  {{\bf E}|X_1|^p}{\sigma ^p }+ 1 \right)
\end{equation}
Let $Z \sim N(0,1)$ and sequence $(X_n)$ be independent. 
We fix $\varepsilon >0$ and choose a compact set $K$, $K \in {\mathcal B}({\bf R})$, such that
\begin{equation}
\begin{array}{lll}
\gamma (K^c) < \varepsilon, &{\bf E}|Z|^p {\bf 1}_{\{Z \in K^c\}} < \varepsilon,& \gamma ( \partial K)=0.
\end{array}
\end{equation}
Hence,
\begin{equation}
\begin{aligned}
\left|{\bf E}|Y_n|^p -{\bf E}|Z|^p \right|& \leq \left|{\bf E}|Y_n|^p {\bf 1}_{\{Z \in K\}} -{\bf E}|Z|^p{\bf 1}_{\{Z \in K\}}\right|+{\bf E}|Y_n|^p {\bf 1}_{\{Z \in K^c\}}+{\bf E}|Z|^p{\bf 1}_{\{Z \in K^c\}}\\
                              & \leq \varepsilon + \varepsilon K(p) \left(\frac  {{\bf E}|X_1|^p}{\sigma ^p }+ 1 \right)+ \varepsilon,\\
			      				\end{aligned}
\end{equation}
for sufficiently large $n$, thanks to (29), (30) and to the convergence of ${\bf E}|Y_n|^p {\bf 1}_{\{Z \in K\}}$ to
${\bf E}|Z|^p{\bf 1}_{\{Z \in K\}}$. Therefore, ${\bf E}|Y_n|^p \to {\bf E}|Z|^p$. This, in particular,
implies the uniform integrability of $(|Y_n|^p)$.
 
Next, we show that for a cost function with $C(x)=O(x^p)$, at infinity, 
all the conditions of Theorem 2 are satisfied. We have ${\bf E}C(2|Z|)< +\infty$,
while the finiteness of $\int C(2|x|) \mu _n (dx) $ follows from
(29) and inequality
\begin{equation}
C(2|Y_n|)=C(2|Y_n|){\bf 1}_{\{|Y_n| \leq x_0\}}+C(2|Y_n|){\bf 1}_{\{|Y_n| > x_0\}}\leq C(x_0){\bf 1}_{\{|Y_n| \leq x_0\}}+ \beta |Y_n|^p{\bf 1}_{\{|Y_n| > x_0\}},
\end{equation}
with a positive constant $ \beta $ such that $C(x) \leq \beta x^p$, for all $x >x_0$.
The CLT provides the weak convergence
$\mu _n \Longrightarrow \gamma $, while we obtain ${\bf E} C(2|Y_n|) \to {\bf E}C(2|Z|)$ 
from the uniform integrability of $(C(2|Y_n|)$ which follows from (32) and uniform integrability
of $(|Y_n|^p)$. Thus, applying Theorem 2, we obtain that $T_c( \mu_n, \gamma ) \to 0$.

\medskip
\noindent
{\bf Remark 2.} Since ${\bf E}|Y_n|^p \to {\bf E}|Z|^p$, the convergence
of ${\bf E} C(|Y_n|)$ to ${\bf E}C(|Z|)$ follows from part (a) and part (c) of the
result of Bickel and Freedman \cite{BF} cited above. Then, ${\bf E} C(2|Y_n|) \to {\bf E}C(2|Z|)$
is implied by the doubling condition (6).

\medskip
\noindent
(ii) Once again, we check that the conditions of Theorem 2 are satisfied.
First,  ${\bf E}C (4 \sqrt 2|Z|) < \infty$ implies that ${\bf E}C (2|Z|) < \infty$ and
that $C(2x)=o(e^{x^2/16})$.

The function $f(x)= \exp ( x^2 /16 )$ has the  expansion
\begin{equation}
f(x)=1+ \sum_{k=1}^{\infty} \frac {(2k-1)!!}{2^{3k}(2k)!}x^{2k}.
\end{equation}
Then Stirling's formula yields, for a generic term $f_k$ of series (33), starting at some index $k_0$
\begin{equation}
f_k=\frac {e^{k+1}(k- \frac 1 2)^{2k-\frac 52}}{2^{4k+1}(k-2)^{k-\frac 32}k^{2k+ \frac 12}}x^{2k}.
\end{equation}
Next,   the Rosenthal inequality (7) can be written, with constants, in the following form (Petrov \cite{Petrov}, inequality (2.35)):
\begin{equation}
{\bf E}|S_n|^k \leq k^kn {\bf E}|X_1|^k+\frac {4k^{\frac k 2+1}}{2^k}e^kn^{\frac k2} \sigma ^k.
\end{equation}
Together (33)--(35) imply, for $Y_n= S_n /(\sigma \sqrt n)$,
\begin{equation}
\begin{aligned}
{\bf E} f(|Y_n|) &\leq M+ \sum _{k=k_0}^{\infty} \frac {e^{k+1}(k-\frac 12)^{2k-\frac 52}{\bf E}|X_1|^{2k}}
{2^{2k+1} \sqrt k (k-2)^{k-\frac 32} \sigma ^{2k}n^{k-1}}+
\sum _{k=k_0}^{\infty} \frac {e^{3k+1}(k-\frac 12)^{2k-\frac 52}}
{2^{5k} \sqrt k (k-2)^{k-\frac 32} k^{k-\frac 12}}\\
&\leq M+ \beta _1 \sum _{k=k_0}^{\infty}k^2{\bf E}|X_1|^{2k}+ \beta_2\sum _{k=k_0}^{\infty} \left( \frac {e^3}{2^5} \right)^k\\
&=M+ \beta _1 Q_1+ \beta _2 Q_2< +\infty
\end{aligned}
\end{equation}
for some positive constants $M$, $\beta _1$, and $\beta _2$.

Since $C(2x)=o(f(x))$, there exists $x_0>0$ such that $C(2x) \leq f(x)$, for all $x>x_0$.
The inequality (36) gives
\begin{equation}
\begin{aligned}
{\bf E} C(2|Y_n|)&={\bf E} C(2|Y_n|){\bf 1}_{\{|Y_n| \leq x_0\}}+{\bf E} C(2|Y_n|){\bf 1}_{\{|Y_n| > x_0\}}\\
                 & \leq C(2x_0)+M+ \beta _1 Q_1+ \beta _2 Q_2.
\end{aligned}
\end{equation}
Therefore, $({\bf E} C(2|Y_n|))$ is bounded and, moreover, uniformly bounded.
Next, we check that $({\bf E} C(2|Y_n|){\bf 1}_A)$, $A \in {\mathcal F}$, is uniformly continuous.

Fix $\varepsilon >0$, and choose $x_1$ positive and such that $C(2x)\leq \varepsilon_1 f(x)$, for all $x>x_1$,
with $\varepsilon _1=\frac {\varepsilon}{2(M+ \beta _1 Q_1+ \beta _2 Q_2)}$.
If $P(A)=\frac {\varepsilon}{2C(2x_1)}$, and in complete similarity to (37) we have
\begin{equation}
\sup _n{\bf E} C(2|Y_n|){\bf 1}_A \leq C(2x_1)P(A)+ \varepsilon _1 \sup _n{\bf E} f(|Y_n|){\bf 1}_A< \varepsilon.
\end{equation}
The inequalities (37) and (38) yield the uniform integrability of $(C(2|Y_n|))$, while
the classical CLT provides $\mu _n\Longrightarrow \gamma $.
Hence, ${\bf E}C(2|Y_n|)\rightarrow {\bf E}C(2|Z|)$ and all the conditions of Theorem 2 are
satisfied. Then $T_c(\mu _n, \gamma) \to 0$.
This concludes the proof of Theorem 3.
\end{proof}

\smallskip
\noindent
{\bf Proof of Theorems 4 and 5. } They are carried out by using the same arguments
as in the proof of Theorem 3.

For sequences of dependent random variables we assume additional conditions of asymptotic independence, 
which yield the moment inequalities (9) and (11) for
the moments  of the sums: the condition (8) for strongly mixing
sequences and the condition (10) for associated sequences.  We also use the bounds on $K(p)$
in (9) and (11)
derived by Doukhan and Louhichi in
\cite {DL}.

\end{document}